\documentclass[11pt]{amsart}
\usepackage{graphicx}
\usepackage[ansinew]{inputenc}    

\newtheorem{theorem}{Theorem}
\theoremstyle{plain}

\newtheorem{corollary}{Corollary}

\newtheorem{definition}{Definition}

\newtheorem{lemma}{Lemma}

\newtheorem{proposition}{Proposition}

\numberwithin{equation}{section}

\begin{document}
\title{Curves with constant curvature ratios}
\author{J. Monterde}
\address{Dep. de Geometria i
Topologia, Universitat de Val\`encia, Avd. Vicent Andr\'es
Estell\'es, 1, E-46100-Burjassot (Val\`encia), Spain}
\email{monterde{@}uv.es}

\thanks{This work was partially supported by a Spanish MCyT grant
BFM2002-00770.}
\date{December 14, 2004}

\subjclass[2000]{Primary 53A04; Secondary 53C40, 53A05}

\keywords{Spherical curves, generalized helices, theorem of
Lancret}

\begin{abstract}
Curves in ${\mathbb R}^n$ for which the ratios between two
consecutive curvatures are constant are characterized by the fact
that their tangent indicatrix is a geodesic in a flat torus. For
$n= 3,4$, spherical curves of this kind are also studied and
compared with intrinsic helices in the sphere.
\end{abstract}

\maketitle

\section{Introduction}

The notion of a generalized helix in ${\mathbb R}^3$, a curve
making a constant angle with a fixed direction, can be generalized
to higher dimensions in many ways. In \cite{RS} the same
definition is proposed but in ${\mathbb R}^n$. In \cite{Ha} the
definition is more restrictive: the fixed direction makes a
constant angle with all the vectors of the Frenet frame. It is
easy to check that this definition only works in the odd
dimensional case. Moreover, in the same reference, it is proven
that the definition is equivalent to the fact that the ratios
$\frac{k_2}{k_1},\frac{k_4}{k_3}, \dots $, $k_i$ being the
curvatures, are constant. This statement is related with the
Lancret Theorem for generalized helices in ${\mathbb R}^3$ (the
ratio of torsion to curvature is constant). Finally, in \cite{Ba}
the author proposes a definition of a general helix in a
$3$-dimensional real-space-form substituting the fixed direction
in the usual definition of generalized helix by a Killing vector
field along the curve.

In this paper we study the curves in ${\mathbb R}^n$ for which all
the ratios $\frac{k_2}{k_1},\frac{k_3}{k_2},\frac{k_4}{k_3}, \dots
$ are constant. We call them curves with constant curvature ratios
or ccr-curves. The main result is that, in the even dimensional
case, a curve has constant curvature ratios if and only if its
tangent indicatrix is a geodesic in the flat torus. In the odd
case, a constant must be added as the new coordinate function.

 In the last section we show
that a ccr-curve in $S^3$ is a general helix in the sense of
\cite{Ba} if and only if it has constant curvatures. To achieve
this result, we have obtained the characterization of spherical
curves in ${\mathbb R}^4$ in terms of the curvatures. Moreover, we
 have also found explicit examples of spherical ccr-curves with
 non-constant curvatures.

\section{Frenet's elements for a curve in ${\mathbb R}^n$}

Let us recall from \cite{Kl} the definition of the Frenet frame
and curvatures.

For $C^{n-1}$ curves, $\alpha$, which have linearly independent
derivatives up to order $n-1$, the moving Frenet frame is
constructed as it were in usual space using the Gram-Schmidt
process. Orthonormal vectors $\{\overrightarrow{{\bf
e_1}},\overrightarrow{{\bf e_2}},\dots, \overrightarrow{{\bf
e_{n-1}}}\}$ are obtained and the last vector is added as the unit
vector in ${\mathbb R}^n$ such that $\{\overrightarrow{{\bf
e_1}},\overrightarrow{{\bf e_2}},\dots, \overrightarrow{{\bf
e_{n}}}\}$ is an orthonormal basis with positive orientation.

The $i$th curvature is defined as
$$k_i= \frac{<\dot{\overrightarrow{{\bf e_i}}},\overrightarrow{{\bf
e_{i+1}}}>}{||\alpha'||},$$ for $i = 1 ,\dots, n-1$.

Frenet's formulae in $n$-space can be written as
\begin{equation}\label{formules-Frenet}
\begin{pmatrix} \dot{\overrightarrow{{\bf e_1}}}(s)\\ \dot{\overrightarrow{{\bf e_2}}}(s)\\
 \dot{\overrightarrow{{\bf e_3}}}(s)\\ \vdots \\ \dot{\overrightarrow{{\bf e_{n-1}}}}(s) \\
 \dot{\overrightarrow{{\bf e_{n}}}}(s) \end{pmatrix} =
\begin{pmatrix}
0 & k_1 & 0 & 0 & \dots & 0 & 0 \\
-k_1 & 0 & k_2& 0 & \dots & 0 & 0 \\
0 & -k_2 & 0 & k_3 & \dots & 0 & 0 \\
\vdots&\vdots&\vdots&\vdots&\ddots&\vdots&\vdots\\
0 & 0 & 0 & 0 & \dots & 0 & k_{n-1} \\
0 & 0 & 0 & 0 & \dots & -k_{n-1} & 0 \\
\end{pmatrix}
\begin{pmatrix} {\overrightarrow{{\bf e_1}}}(s)\\ {\overrightarrow{{\bf e_2}}}(s)\\
 {\overrightarrow{{\bf e_3}}}(s)\\ \vdots \\ {\overrightarrow{{\bf e_{n-1}}}}(s) \\
 {\overrightarrow{{\bf e_{n}}}}(s) \end{pmatrix}.
 \end{equation}

In accordance with \cite{RS} we will say that a curve is twisted
if its last curvature, $k_{n-1}$ is not zero. Sometimes, we will
also say that the curve is not regular.

\section{ccr-curves}

Instead of looking for curves making a constant angle with a fixed
direction as in \cite{Ha} or \cite{RS}, we will study another way
of generalizing the notion of helix.

\begin{definition}
A curve $\alpha:I\to {\mathbb R}^n$ is said to have constant
curvature ratios (that is to say, it is a ccr-curve) if all the
quotients $\frac{k_{i+1}}{k_i}$ are constant.
\end{definition}

As is well known, generalized helices in ${\mathbb R}^3$ are
characterized by the fact that the quotient $\frac\tau\kappa$ is
constant (Lancret's theorem). It is in this sense that ccr-curves
are a generalization to ${\mathbb R}^n$ of generalized helices in
${\mathbb R}^3$.

In \cite{Ha} the author defines a generalized helix in the
n-dimensional space (n odd) as a curve satisfying that the ratios
$\frac{k_2}{k_1},\frac{k_4}{k_3}, \dots $ are constant. It is also
proven that a curve is a generalized helix if and only if there
exists a fixed direction which makes constant angles with all the
vectors of the Frenet frame. Obviously, ccr-curves are a subset of
generalized helices in the sense of \cite{Ha}.

\subsection{Examples}

\subsubsection{Example with constant curvatures}
The subset of ${\mathbb R}^{2n}$ parametrized by
$\overrightarrow{\bf x}(u_1,u_2,\dots,u_n) = $
$$= (r_1\cos(u_1),r_1\sin(u_1),
r_2\cos(u_2),r_2\sin(u_2),\dots,r_n\cos(u_n),r_n\sin(u_n))$$ where
$u_i\in{\mathbb R}$ is called a flat torus in ${\mathbb R}^{2n}$.

By analogy, the subset of ${\mathbb R}^{2n+1}$ parametrized by
$\overrightarrow{\bf x}(u_1,u_2,\dots,u_n) = $
$$= (r_1\cos(u_1),r_1\sin(u_1),
r_2\cos(u_2),r_2\sin(u_2),\dots,r_n\cos(u_n),r_n\sin(u_n),a)$$
where $u_i\in{\mathbb R}$ and $a$ is a real constant, will be
called a flat torus in ${\mathbb R}^{2n+1}$.

It is just a matter of computation to show that any curve in a
flat torus of the kind $$\alpha(t) = \overrightarrow{\bf x}(m_1
t,m_2 t,\dots,m_n t)$$ has all its curvatures constant (see
\cite{Ro}).

These curves are the geodesics of the flat tori, and it is proven
in the cited paper that they are twisted curves if and only if the
constants $m_i\ne m_j$ for all $i\ne j$.

\subsubsection{Example with non-constant curvatures}

Now, let $k(s)$ be a positive function. Let us define $g(s) =
\int_0^sk(u) du$. If $\alpha$ is a curve parametrized by its
arc-length and with constant curvatures, $a_1,a_2,\dots,a_{n-1}$,
then the curve $\beta(s) = \int_0^s \overrightarrow{\bf
e_1}^\alpha(g(u)) du$ is a curve whose curvatures are $k_i(s) =
a_i k(s)$.

Note that $\dot{\beta}(s) = \overrightarrow{\bf
e_1}^\alpha(g(s))$. This implies that $\overrightarrow{\bf
e_1}^\beta(s) = \overrightarrow{\bf e_1}^\alpha(g(s))$. Taking
derivatives $k_1^\beta(s) \overrightarrow{\bf e_2}^\beta(s) =
k_1^\alpha(g(s)) \overrightarrow{\bf e_2}^\alpha(g(s))k(s)$.
Therefore,
$$\overrightarrow{\bf
e_2}^\beta(s) = \overrightarrow{\bf e_2}^\alpha(g(s)), \qquad
\text{and}\qquad k_1^\beta(s) = a_1 k(s).$$

By similar arguments it is possible to show that $k_i^\beta(s) =
a_i k(s)$ for any $i =1 ,\dots, n-1$. Therefore, $\beta$ is a
ccr-curve with non-constant curvatures.

In the next section we will show that any ccr-curve is of this
kind.

\section{Solving the natural equations for ccr-curves}

The Frenet formulae can be explicitly integrated only for some
particular cases. Ccr-curves are one of these. In fact, Frenet's
formulae are

$$
\begin{pmatrix} \dot{\overrightarrow{{\bf e_1}}}(s)\\ \dot{\overrightarrow{{\bf e_2}}}(s)\\
 \dot{\overrightarrow{{\bf e_3}}}(s)\\ \vdots \\ \dot{\overrightarrow{{\bf e_{n-1}}}}(s) \\
 \dot{\overrightarrow{{\bf e_{n}}}}(s) \end{pmatrix} = k_1(s)
\begin{pmatrix}
0 & 1 & 0 & 0 & \dots & 0 & 0 \\
1 & 0 & c_2& 0 & \dots & 0 & 0 \\
0 & -c_2 & 0 & c_3 & \dots & 0 & 0 \\
\vdots&\vdots&\vdots&\vdots&\ddots&\vdots&\vdots\\
0 & 0 & 0 & 0 & \dots & 0 & c_{n-1} \\
0 & 0 & 0 & 0 & \dots & -c_{n-1} & 0 \\
\end{pmatrix}
\begin{pmatrix} {\overrightarrow{{\bf e_1}}}(s)\\ {\overrightarrow{{\bf e_2}}}(s)\\
 {\overrightarrow{{\bf e_3}}}(s)\\ \vdots \\ {\overrightarrow{{\bf e_{n-1}}}}(s) \\
 {\overrightarrow{{\bf e_{n}}}}(s) \end{pmatrix},$$
 for some constants, $c_2,\dots, c_{n-1}$.

Reparametrization of the curve allows that system to be reduced to
an easier one. The reparametrization is given by the inverse
function of
$$g(s) = \int_0^s k_1(u)du.$$
Note that $t=g(s)$ is a reparametrization because $k_1$ is a
positive function. The reparametrization we need is the inverse
function $s = g^{-1}(t)$. It is a simple matter to verify that,
with respect to parameter $t$, the Frenet's formulae are reduced
to a linear system of first order differential equations with
constant coefficients

\begin{equation}\label{eqs-coefs-const}
\begin{pmatrix} \overrightarrow{\bf e_1}'(t)\\ \overrightarrow{\bf e_2}'(t)\\
\overrightarrow{\bf e_3}'(t)\\ \vdots \\ \overrightarrow{\bf e_{n-1}}'(t) \\
\overrightarrow{\bf e_n}'(t) \end{pmatrix} =
\begin{pmatrix}
0 & 1 & 0 & 0 & \dots & 0 & 0 \\
1 & 0 & c_2& 0 & \dots & 0 & 0 \\
0 & -c_2 & 0 & c_3 & \dots & 0 & 0 \\
\vdots&\vdots&\vdots&\vdots&\ddots&\vdots&\vdots\\
0 & 0 & 0 & 0 & \dots & 0 & c_{n-1} \\
0 & 0 & 0 & 0 & \dots & -c_{n-1} & 0 \\
\end{pmatrix}
\begin{pmatrix} {\overrightarrow{{\bf e_1}}}(t)\\ {\overrightarrow{{\bf e_2}}}(t)\\
 {\overrightarrow{{\bf e_3}}}(t)\\ \vdots \\ {\overrightarrow{{\bf e_{n-1}}}}(t) \\
 {\overrightarrow{{\bf e_{n}}}}(t) \end{pmatrix}.
\end{equation}

We can apply the well-known methods of integration of systems of
linear equations with constant coefficients. Let $F_n$ be the
matrix of constant coefficients of this system.

\subsection{Eigenvalues and their multiplicity} The first thing we have to do is to compute the
eigenvalues of the coefficient matrix.

Due to the skewsymmetry of the matrix, it can have not real
eigenvalues other than zero. Due to the fact that the determinant
of $F_n$ vanishes only for odd $n$, we can say that for odd
dimensions, $0$ is an eigenvalue, whereas for even dimensions, $0$
is an eigenvalue only if $k_{n-1} = 0$.

By definition, we have that constants $c_2,c_3,\dots,c_{n-2}$ are
not zero. If the last constant, $c_{n-1}$, vanishes, then the same
happens with the last curvature function $k_{n-1}$. In this case
the curve is included in a hyperspace, so we can consider it to be
a curve in an $n-1$ dimensional space.

Therefore, from now on, we shall consider that all the curvatures,
and then all the constants $c_i$, are not zero.

Note that, in this case, for any $x\in {\mathbb C}$, the rank (in
${\mathbb C}$) of the matrix
$$\begin{pmatrix}
x & 1 & 0 & 0 & \dots & 0 & 0 \\
1 & x & c_2& 0 & \dots & 0 & 0 \\
0 & -c_2 & x & c_3 & \dots & 0 & 0 \\
\vdots&\vdots&\vdots&\vdots&\ddots&\vdots&\vdots\\
0 & 0 & 0 & 0 & \dots & x & c_{n-1} \\
0 & 0 & 0 & 0 & \dots & -c_{n-1} & x \\
\end{pmatrix}$$
is at least $n-1$. Therefore, their eigenvalues are all of
multiplicity $1$.

\subsection{Canonical Jordan form}

Let $a_\ell\pm {\bf i} b_\ell$, $\ell = 1, \dots, [\frac n2]$,
with $a_\ell,b_\ell\in{\mathbb R}$, be the non-zero eigenvalues of
the coefficient matrix. Therefore, for $n= 2k$, the associated
canonical Jordan form is of the kind
$$\begin{pmatrix}
J_1 & 0  & \dots  & 0 \\
0 & J_2  & \dots  & 0 \\
\vdots & \vdots & \ddots & \vdots\\
0 & 0 & \dots  &J_k \\
\end{pmatrix}$$
where $J_\ell = \begin{pmatrix}  a_\ell & -b_\ell \\
b_\ell & a_\ell \end{pmatrix}$.

The matrix can be diagonalized because all the eigenvalues are of
multiplicity one. Therefore, there is a orthogonal matrix, $S$,
such that if $C$ is the matrix of constant coefficients, then
$$C = S^{-1}JS.$$

 Therefore, the general solution of
the system for the first vector is
$${\overrightarrow{\bf e}_1}(u) :=
\sum_{\ell=1}^{k}\overrightarrow{A_\ell}\ e^{a_\ell u}\cos(b_\ell\
u) +\overrightarrow{B_\ell}\ e^{a_\ell u} \sin(b_\ell\ u),$$ where
$\{\overrightarrow{A_\ell},\overrightarrow{B_\ell}\}_{\ell = 1}^k$
is a family of orthogonal vectors.

For $n= 2k+1$, the associated canonical Jordan form is of the kind
$$\begin{pmatrix}
0 & 0 & 0  & \dots  & 0 \\
0 & J_1 & 0  & \dots  & 0 \\
0 & 0 & J_2  & \dots  & 0 \\
\vdots & \vdots & \vdots & \ddots & \vdots\\
0 & 0 & 0 & \dots  &J_k \\
\end{pmatrix}$$

Now, the general solution of the system for the first vector is
$${\overrightarrow{\bf e}_1}(u) := \overrightarrow{A_0}+
\sum_{\ell=1}^{k}\overrightarrow{A_\ell}\ e^{a_\ell u}\cos(b_\ell\
u) +\overrightarrow{B_\ell}\ e^{a_\ell u}\sin(b_\ell\ u),$$ where
$\{\overrightarrow{A_0}\}\cup
\{\overrightarrow{A_\ell},\overrightarrow{B_\ell}\}_{\ell = 1}^k$
is a family of orthogonal vectors.

\subsection{The eigenvalues are pure imaginaries}

Condition $||{\overrightarrow{\bf e}_1}(u)|| = 1$ for all $u$
implies that all the real parts of the eigenvalues are zero.
Indeed, if, for example, $a_1 \ne 0$, then let $m$ be a non-zero
coordinate of $\overrightarrow{A_1}$.

Bearing in mind that
$$|m|\ e^{a_1 u}\ |\cos(b_1 u)| \le ||{\overrightarrow{\bf
e}_1}(u)||,$$ and that the left-hand member is an unbounded
function, then $||{\overrightarrow{\bf e}_1}(u)||\ne 1$.

Therefore, all the real parts of the eigenvalues are zero and the
general solution (in the even case) of the system for the first
vector is
$${\overrightarrow{\bf e}_1}(u) :=
\sum_{\ell=1}^{k}\overrightarrow{A_\ell} \cos(b_\ell\ u)
+\overrightarrow{B_\ell} \sin(b_\ell\ u).$$

Analogously for the odd case.

Moreover, let us recall that the vectors
$\{\overrightarrow{A_i},\overrightarrow{B_i}\}_{i=1}^k$ are an
orthogonal base of ${\Bbb R}^n$ associated to the canonical Jordan
form.

\subsection{The main result}

Finally, an isometry of ${\mathbb R}^n$ allows us to state the
next result.

\begin{theorem}  A curve has constant curvature ratios if and only
if its tangent indicatrix is a twisted geodesic on a flat torus.
\end{theorem}

Note that in the odd dimensional case this result implies that the
last coordinate of the tangent indicatrix is a constant. So, there
is a direction making a constant angle with the curve.
Nevertheless, this is not the case in the even dimensional case.
There are no fixed directions making a constant angle with the
tangent vector.

When all the curvatures are constant, then the curve is also a
ccr-curve and its tangent indicatrix is of the kind described in
the previous statement. Moreover, the reparametrization $g(s) =
\int_0^s k_1(u)du$ is just the product by a constant.

Since the integration of a geodesic on a flat torus in ${\mathbb
R}^{2k}$ with respect to its parameter is again a curve of the
same kind, we get the following corollary:

\begin{corollary}  A curve has constant curvatures if and only if
it is \begin{enumerate} \item a twisted geodesic on a flat torus,
in the even dimensional case, or \item a twisted geodesic on a
flat torus times a linear function of the parameter, in the odd
dimensional case.
\end{enumerate}
\end{corollary}

\subsection{$n = 3$}

The eigenvalues of the matrix of coefficients are $0$ and $\pm
\sqrt{1+c^2}\ {\bf i}$  ($c= c_2$, to simplify).

Therefore, the general solution of the system for the first vector
is
$${\overrightarrow{{\bf e_1}}}(u) =\overrightarrow{A_1}
+\overrightarrow{A_2} \cos(\sqrt{1+c^2} u) +\overrightarrow{A_3}
\sin(\sqrt{1+c^2} u),$$ where $\overrightarrow{A_i}, i = 1,2,3$
are constant vectors.

Once we have the tangent vector, we only have to undo the
reparametrization and to integrate to obtain the curve

$$\alpha(s) = x_0+ \overrightarrow{c_1} s +\overrightarrow{c_2}\int_0^s \cos(\sqrt{1+c^2}g(v))dv +
\overrightarrow{c_3} \int_0^s\sin(\sqrt{1+c^2}g(v))dv.$$

\subsection{$n = 4$}

The eigenvalues are
$$\pm \frac {\bf i}{\sqrt{2}} \sqrt{ 1+ c_2^2 + c_3^2 \pm \sqrt{(1 + c_2^2 +
c_3^2)^2-4 c_3^2}}.$$

Therefore, the general solution of the system for the first vector
is
$${\overrightarrow{\bf e}_1}(u) :=
\overrightarrow{A_1} \cos(m_{+} u) +\overrightarrow{B_1}
\sin(m_{+} u)+\overrightarrow{A_2} \cos(m_{-} u)
+\overrightarrow{B_2} \sin(m_{-} u),$$ where
$$
m_{\pm} = \frac 1{\sqrt{2}} \sqrt{ 1+ c_2^2 + c_3^2 \pm \sqrt{(1 +
c_2^2 + c_3^2)^2-4 c_3^2}}
$$
and where $\overrightarrow{A_i},\overrightarrow{B_i}, i = 1,2$ are
constant vectors.

\section{Spherical ccr-curves}

In order to compare ccr-curves with the definition of generalized
helices given in \cite{Ba}, we will try to determine which
ccr-curves are included in a sphere.

\begin{lemma}
A curve $\alpha:I\to {\mathbb R}^4$ is spherical, i.e., it is
contained in a sphere of radius $R$, if and only if
\begin{equation}\label{spherical-eq-n=4}
\frac 1{k_1^2} + \left(\frac{\dot{k_1}}{k_1^2k_2} \right)^2 +
\frac
1{k_3^2}\left(\left(\frac{\dot{k_1}}{k_1^2k_2}\right)^{\dot{\
}}-\frac {k_2}{k_1} \right)^2 = R^2.
\end{equation}
\end{lemma}

\begin{proof}
The proof here is similar to that for spherical curves in
${\mathbb R}^3$. It consists in obtaining information thanks to
successive derivatives of the expression
$<\alpha(s)-m,\alpha(s)-m> = R^2$, where $m$ is the center of the
sphere. In particular, what can be proven is that spherical curves
can be decomposed as
\begin{equation}\label{spherical-curve}
\alpha(s) = m - \frac R{k_1} {\overrightarrow{{\bf e_2}}}(s) +R
\frac{\dot{k_1}}{k_1^2k_2} {\overrightarrow{{\bf e_3}}}(s) +R\frac
1{k_3}\left(\left(\frac{\dot{k_1}}{k_1^2k_2}\right)^\cdot +\frac
{k_2}{k_1} \right){\overrightarrow{{\bf e_4}}}(s).
\end{equation}
\end{proof}

As a corollary we obtain the classical result for spherical
three-dimensional curves:

\begin{corollary}
A curve $\alpha:I\to {\mathbb R}^3$ is spherical, i.e., it is
contained in a sphere of radius $R$, if and only if
\begin{equation}\label{spherical-eq-n=3}
\frac 1{k_1^2} + \left(\frac{\dot{k_1}}{k_1^2k_2} \right)^2 = R^2.
\end{equation}
\end{corollary}

From now on, we shall suppose that $m= 0$ and $R= 1$.

\subsection{Spherical ccr-curves in ${\mathbb R}^3$}

In this case, we can rewrite Eq. \ref{spherical-eq-n=3} in terms
of curvature, $k_1 = \kappa$, and torsion $k_2 = \tau = c \kappa$,
$c$ being a constant.

$$\frac {\dot\kappa}{\kappa^2\sqrt{\kappa^2-1}} = \pm c.$$

Let us consider just the positive sign. This differential equation
can be integrated and the solution is
$$\kappa(s) = \frac 1{\sqrt{1-(cs+s_0)^2}}.$$

Thanks to a shift of the parameter we get that the curvature and
torsion of a spherical generalized helix are given by
$$\kappa(s) = \frac 1{\sqrt{1-c^2s^2}},\qquad
\tau(s) = \frac c{\sqrt{1-c^2s^2}}.$$

We now need to compute the reparametrization
$$u=g(s) = \int_{0}^s \kappa(t)dt = \frac1c \arcsin(cs).$$

With the appropriate initial conditions, the generalized spherical
helix is
$$\begin{array}{rcl}
\alpha_c(s) &=& (\sqrt{1-c^2 s^2}
\cos(\frac{\sqrt{1+c^2}\arcsin(cs)}{c})+ \frac{c^2 s}{\sqrt{1+c^2}}\sin(\frac{\sqrt{1+c^2}\arcsin(cs)}{c}),\\[2mm]
&&-\sqrt{1-c^2 s^2} \sin(\frac{\sqrt{1+c^2}\arcsin(cs)}{c})+
\frac{c^2 s}{\sqrt{1+c^2}}\cos(\frac{\sqrt{1+c^2}\arcsin(cs)}{c}),\\[2mm]
&& \frac{c s}{\sqrt{1+c^2}}) \end{array}
$$

Note that the curve $\alpha_c$ is defined in the interval $]-\frac
1c,\frac 1c[$. If we change the parameter in accordance with $s =
\frac 1c \sin t$, the spherical helix is now parametrized as
$$\begin{array}{rcl}
\beta_c(t) &=& (\cos t
\cos(\frac{\sqrt{1+c^2}}{c}t)+ \frac{c}{\sqrt{1+c^2}}\sin t\sin(\frac{\sqrt{1+c^2}}{c}t),\\[2mm]
&&-\cos t \sin(\frac{\sqrt{1+c^2}}{c}t)+
\frac{c}{\sqrt{1+c^2}}\sin
t\cos(\frac{\sqrt{1+c^2}}{c}t),\frac{\sin t}{\sqrt{1+c^2}})
\end{array}
$$

Now, it is clear that the projection of these curves on the plane
$xy$ are arcs of epicycloids. This result was known by W.
Blaschke, as is mentioned in \cite{St}, where it is also proven by
different methods.

\subsection{Spherical ccr-curves in ${\mathbb R}^4$}

\subsubsection{The constant curvatures case.}

The curve
$$\alpha(s) = \frac 1{\sqrt{r_1^2+r_2^2}} (\frac{r_1}{m_1}
\sin(m_1 s),-\frac{r_1}{m_1} \cos(m_1 s),\frac{r_2}{m_2} \sin(m_2
s),-\frac{r_2}{m_2} \cos(m_2 s))$$ is a spherical curve (with
radius $1$), if and only if
$$r_1^2m_2^2+r_2^2m_1^2 = m_1^2m_2^2(r_1^2+r_2^2).$$
\subsubsection{The non-constant case.}

In this case, we can rewrite Eq. \ref{spherical-eq-n=4} in terms
of curvature, $k_1$, $k_2 = c_2 k_1$ and $k_3 = c_3 k_1$, where
$c_2, c_3$  are constants.
\begin{equation}
\frac 1{k_1^2} + \left(\frac{\dot{k_1}}{c_2k_1^3} \right)^2 +
\frac
1{c_3^2k_1^2}\left(\left(\frac{\dot{k_1}}{c_2k_1^3}\right)^\cdot
+c_2 \right)^2 = 1.
\end{equation}
By changing $f = \frac 1{k_1^2}$ the equation is reduced to
\begin{equation}\label{reduced}  f + \frac 1{4c_2^2} \dot{f}^2 +\frac 1{c_3^2} f(-\frac 1{2c_2}
\ddot{f} + c_2)^2 = 1.\end{equation}
 Computation of the general
solution seems to be a difficult task. Instead, we can try to
compute some particular solutions.

For instance, the constant solution $f(s) = \frac{c_3^2}{c_2^2 +
c_3^2}$ or the polynomial solutions of degree $2$
$$f(s) =
\frac{-2 c_2^2 + c_3^2 - c_3\sqrt{-8 c_2^2 + c_3^2}}{2(c_2^2 +
c_3^2)} + \frac 12
   \left(2 c_2^2 - c_3^2 - c_3\sqrt{-8 c_2^2 + c_3^2}\right) s^2,$$
$$f(s) = 2 c_2 s +\frac 12
   \left(2 c_2^2 - c_3^2 - c_3\sqrt{-8 c_2^2 + c_3^2}\right)
   s^2.$$
For these three particular solutions the reparametrization $g$,
where $g(s) = \int_0^s k_1(t)dt = \int_0^s \frac
1{\sqrt{f(t)}}dt,$ can be computed explicitly. We can thus obtain
explicit examples of ccr-curves in $S^3$ with non-constant
curvatures.

A particular case. With $c_2= \frac 12, c_3:=\frac{\sqrt{3}}2$,
then $m_1= \sqrt{\frac 32}, m_2 = \frac 1{\sqrt{2}}$ and $r_1= r_2
= \frac 1{\sqrt{2}}$. The function $f(s) = \frac 12 - 2s^2$ is a
solution of Eq. \ref{reduced}. Therefore, $k_1(s) = \frac
2{\sqrt{1-4s^2}}$, and $g(s) = \int_0^s \frac 2{\sqrt{1-4t^2}}dt =
\arcsin(2s)$.

If $${\overrightarrow{{\bf e_1}}}(t) = \frac 1{\sqrt{2}}
(\cos(\sqrt{\frac 32}t), \sin(\sqrt{\frac 32}t),\cos(\frac
1{\sqrt{2}}t),\sin(\frac 1{\sqrt{2}}t)),$$ then
$$\alpha(s) = (0, -\frac{\sqrt{3}}2, 0, \frac 12)+\int_0^s
{\overrightarrow{{\bf e_1}}}(\arcsin(2u))du, \quad s\in\ ]-\frac
12,\frac 12[$$ is a spherical ccr-curve with center at the origin
of coordinates, with radius $1$ and with non-constant curvatures.

\section{Intrinsic generalized helices}

In \cite{Ba} the author proposes a definition of general helix on
a $3$-dimensional real-space-form substituting the fixed direction
in the usual definition of generalized helix by a Killing vector
field along the curve.

Let $\alpha:I\to M$ be an immersed curve in a $3$-dimensional
real-space-form $M$. Let us denote  the intrinsic Frenet frame by
$\{\overrightarrow{{\bf t}}, \overrightarrow{{\bf
n}},\overrightarrow{{\bf b}}\}$. The intrinsic Frenet's formulae
are
\begin{equation}\label{intrinsic-Frenet-formulae}
\left\{\aligned
\nabla_{\overrightarrow{{\bf t}}}\overrightarrow{{\bf t}} &= \kappa \overrightarrow{{\bf n}},\\
\nabla_{\overrightarrow{{\bf t}}}\overrightarrow{{\bf n}} &=
-\kappa \overrightarrow{{\bf t}}+
\tau\overrightarrow{{\bf b}},\\
\nabla_{\overrightarrow{{\bf t}}}\overrightarrow{{\bf b}} &= -\tau
\overrightarrow{{\bf n}},
\endaligned
\right.
\end{equation}
where $\nabla$ is the Levi-Civita connection of $M$ and where
$\kappa$ and $\tau$ are called the intrinsic curvature and torsion
functions of curve $\alpha$, respectively.

From now on we shall suppose that $M= S^3$. Therefore, any curve
on $S^3$ can also be considered to be a curve in ${\mathbb R}^4$.
We shall try to obtain the relationship between the Frenet
elements, $\{\overrightarrow{{\bf e_1}},\overrightarrow{{\bf
e_2}},\overrightarrow{{\bf e_3}},\overrightarrow{{\bf e_4}}, k_1,
k_2,k_3\}$, of the curve as a curve in $4$-dimensional Euclidian
space and the intrinsic Frenet elements $\{\overrightarrow{{\bf
t}}, \overrightarrow{{\bf n}},\overrightarrow{{\bf b}},
\kappa,\tau\}$. Note first that $ \overrightarrow{{\bf t}}=
\overrightarrow{{\bf e_1}}.$ Then
$$ \nabla_{\overrightarrow{{\bf t}}}\overrightarrow{{\bf t}} =
\dot{\overrightarrow{\bf e_1}}- <\dot{\overrightarrow{\bf
e_1}},\alpha>\alpha = k_1(\overrightarrow{\bf e_2}-
<{\overrightarrow{\bf e_2}},\alpha>\alpha),$$ where we have used
 as the Gauss map of the sphere the identity map.

Therefore
\begin{equation}\label{intrinsic-normal} \overrightarrow{\bf n} =
\frac{\nabla_{\overrightarrow{{\bf t}}}\overrightarrow{{\bf
t}}}{||\nabla_{\overrightarrow{{\bf t}}}\overrightarrow{{\bf
t}}||}=\frac 1{\sqrt{1-<{\overrightarrow{\bf e_2}},\alpha>^2}}
(\overrightarrow{\bf e_2}- <{\overrightarrow{\bf
e_2}},\alpha>\alpha),
\end{equation}
 and
 $$\kappa = <\nabla_{\overrightarrow{{\bf t}}}\overrightarrow{{\bf
 t}},\overrightarrow{\bf n}> = k_1 \sqrt{1-<{\overrightarrow{\bf
 e_2}},\alpha>^2}= \sqrt{k_1^2-1},$$
which were obtained using  Eq. \ref{spherical-curve}.

The intrinsic binormal vector is the only vector such that
$\{\overrightarrow{{\bf t}}, \overrightarrow{{\bf
n}},\overrightarrow{{\bf b}}, \alpha\}$ is an orthonormal basis of
${\mathbb R}^4$ with positive orientation. Then
$$ \overrightarrow{{\bf b}} = \alpha\wedge \overrightarrow{{\bf
t}}\wedge \overrightarrow{{\bf n}}.$$ Now, by replacing the
intrinsic tangent and normal with $ \overrightarrow{{\bf t}}=
\overrightarrow{{\bf e_1}}$ and \ref{intrinsic-normal}, we get
$$ \overrightarrow{{\bf b}} = \frac {k_1}{\sqrt{k_1^2-1}}\ \alpha\wedge \overrightarrow{{\bf
e_1}}\wedge \overrightarrow{{\bf e_2}} = \frac
{1}{\sqrt{1-(\frac{1}{k_1})^2}}\ \alpha\wedge \overrightarrow{{\bf
e_1}}\wedge \overrightarrow{{\bf e_2}}.$$ Therefore
$$\dot{\overrightarrow{\bf b}} = \left(\frac
{1}{\sqrt{1-(\frac{1}{k_1})^2}}\right)^\cdot\ \alpha\wedge
\overrightarrow{{\bf e_1}}\wedge \overrightarrow{{\bf e_2}}+ \frac
{1}{\sqrt{1-(\frac{1}{k_1})^2}}\ \alpha\wedge \overrightarrow{{\bf
e_1}}\wedge k_2\overrightarrow{{\bf e_3}}.$$ A consequence of this
computation is that $<\dot{\overrightarrow{\bf b}}, \alpha> = 0$,
and therefore, $\nabla_{\overrightarrow{{\bf
t}}}\overrightarrow{{\bf b}}= \dot{\overrightarrow{\bf b}}$.
Finally,
$$\begin{array}{rcl}
\tau &=& - <\nabla_{\overrightarrow{{\bf t}}}\overrightarrow{{\bf
b}},\overrightarrow{{\bf n}}> \\[3mm]
&=& -<\frac {1}{\sqrt{1-(\frac{1}{k_1})^2}}\ \alpha\wedge
\overrightarrow{{\bf e_1}}\wedge k_2\overrightarrow{{\bf
e_3}},\frac 1{\sqrt{1-(\frac{1}{k_1})^2}} \overrightarrow{\bf
e_2}>\\[5mm]
&=& -\frac {k_2}{1-(\frac{1}{k_1})^2}<\alpha\wedge
\overrightarrow{{\bf e_1}}\wedge \overrightarrow{{\bf e_3}},
\overrightarrow{\bf e_2}>=  \frac {k_2}{1-(\frac{1}{k_1})^2} =
\frac {k_1^2k_2}{\kappa^2}.
\end{array}$$
\begin{proposition}
The only $4$-dimensional spherical non-trivial ccr-curves which
are also intrinsic generalized helices of $S^3$ are helices, i.e.,
curves with all curvatures constant.
\end{proposition}
\begin{proof}
As it is proven in \cite{Ba},  a curve in $S^3$ is an intrinsic
helix if and only if $\tau = 0$ or there exists a constant $b$
such that $\tau = b\kappa\pm 1$.

The case $\tau = 0$ implies that $k_1k_2 = 0$ and we get a
non-regular curve.

In the other case, if the curve is also a ccr-curve (with $k_2 = c
k_1$), then
$$\frac
{ck_1^3}{\kappa^2} = b \kappa \pm 1.$$ Equivalently
$$(\frac
{ck_1^3}{k_1^2-1}\mp 1)^2 = b (k_1^2-1).$$ That is, the function
$k_1$ is the solution of a polynomial equation with constant
coefficients; and, therefore, the function $k_1$ is constant, and
so the other two curvatures $k_2$ and $k_3$ are also constant. The
same happens with $\kappa$ and $\tau$. We are then in the presence
of a helix according to the designation in \cite{Ba}, or a
geodesic in a flat torus in ${\mathbb R}^4$ according to
\cite{Ro}.
\end{proof}


\begin{thebibliography}{99}

\bibitem{Ba} M. Barros, {\it General helices and a theorem of
Lancret}, Proceedings of the Am. Math. Soc., vol. 125, 1503-1509
(1997).

\bibitem{dC}  do Carmo, M. P., {\it Differential Geometry of
curves and surfaces}, Prentice-Hall Inc., (1976).

\bibitem{FGL} A. Ferrández, A. Giménez, P. Lucas, {\it Null
generalized helices in Lorentz-Minkowski spaces}, J. of Physics,
A: Math. and general, vol 35, 8243-8251 (2002).

\bibitem{Ha}  H. A. Hayden, {\it On a generalized helix in a
Riemannian $n$-space}, Proc. London Math. Soc., vol 32, 37-45
(1931).

\bibitem{Kl}  W. Klingenberg, {\it A course in Differential
Geometry}, Springer-Verlag (1978).

\bibitem{Ro} S. Rodrigues Costa, {\it On closed twisted curves},
Proceedings of the Am. Math. Soc., vol. 109, 205-214 (1990).

\bibitem{RS} M.C. Romero-Fuster, E. Sanabria-Codesal, {\it
Generalized helices, twistings and flattenings of curves in
$n$-space}, Matemática Contemporânea, vol. 17, 267-280 (1999).

\bibitem{St} Struick, D. J.; {\it Lectures on Classical
Differential Geometry}, Dover, New-York, (1988).

\end{thebibliography}
\end{document}